\documentclass[12pt]{article}
\usepackage{amsmath,amssymb,amsthm, amscd}

\newtheorem{Prop}{Proposition}
\newtheorem{Lem}[Prop]{Lemma}

\newtheorem{Cor}[Prop]{Corollary}

\newtheorem{Claim}[Prop]{Claim}

\begin{document}

\begin{center}
{\bf A remark on isolated symplectic singularities with trivial local fundamental group}
\end{center}
\vspace{0.4cm}

\begin{center}
{\bf Yoshinori Namikawa}      
\end{center}
\vspace{0.2cm}

\begin{abstract} Recently, Bellamy et al. constructed an infinite series of 4-dimensional isolated symplectic sngularities with trivial local fundamental group, inspired by a question of Beauville. In this short note, we introduce an easy construction of isolated symplectic singularities (of any dimension) with trivial local fundamental group. We will use a toric hyperk\"{a}hler construction. 

MSC2020: 14B05, 14J17, 14L30, 32M05, 53D20, 53C26
\end{abstract}
 
\section{Introduction} 

Recently, Bellamy et al. \cite{BBFJLS} constructed an infinite series of 4-dimensional isolated symplectic singularities 
with trivial local fundamental group, inspired by a question of Beauville \cite{Be}. In this short note, we introduce an easy construction  of  isolated 
symplectic singularities (of any dimension) with trivial local fundamental group.  We will use a toric hyperk\"{a}hler construction (\cite{BD}, \cite{Go}, \cite{HS}, \cite{Ko}).   
 
Let $n$ be a positive integer with $n \geq 3$ and 
let $\mathbf{a} := (a_1, ..., a_n)$ be an $n$-tuple of nonzero integers such that $\mathrm{GCD}(a_i, a_j) = 1$ for any 
$i \ne j$.    
Then $\mathbf{C}^*$ acts on the affine space $\mathbf{C}^{2n}$ by 
$$(z_1, ..., z_n, w_1, ..., w_n) \to (t^{a_1}z_1, ..., t^{a_n}z_n, t^{-a_1}w_1, ..., t^{-a_n}w_n).$$ 
By rearranging the order of $z_1, ..., z_n, w_1, ..., w_n$, we may assume that $a_i$ are all positive and 
$a_1 \le a_2 \le ... \le a_n$. 
Define a symplectic form  $\omega$ on $\mathbf{C}^{2n}$ by $$\omega := \sum_{1 \le i \le n} dw_i \wedge dz_i.$$  
Then the $\mathbf{C}^*$-action on $(\mathbf{C}^{2n}, \omega)$ is a Hamiltonian action and 
its moment map $\mu$ is given by 
$$\mathbf{C}^{2n} \to \mathbf{C}\:\:\: (z_1, ..., z_n, w_1, ..., w_n) \to \sum_{1 \le i \le n}a_iz_iw_i. $$
Define $Y(\mathbf{a}, 0) := \mu^{-1}(0)/\hspace{-0.1cm}/\mathbf{C}^* (= \mathrm{Spec}\: \mathbf{C}[\mu^{-1}(0)]^{\mathbf{C}^*})$ and we call it a toric hyperk\"{a}hler variety associated with 
$\mathbf{a}$.  

In this article we shall prove: \vspace{0.2cm} 

\begin{Prop}\label{1}  The toric hyperk\"{a}hler variety $Y(\mathbf{a}, 0)$ has an isolated symplectic singularity of dimension 
$2n -2$ with $\pi_1(Y(\mathbf{a},0)_{reg}) = 1$. If $\mathbf{a} \ne \mathbf{a}'$, then $Y(\mathbf{a}, 0) \ne 
Y(\mathbf{a}', 0)$. 
\end{Prop}
    
In the following we give a proof of Proposition \ref{1}. 
Let $\pi: \mu^{-1}(0) \to Y(\mathbf{a}, 0)$ be the quotient map.  
Put $\bar{0} := \pi (0)$. \vspace{0.2cm}

\begin{Claim}\label{2} The variety $Y(\mathbf{a}, 0)$ is a symplectic variety of dimension $2n-2$ with an  isolated singularity at $\bar{0}$. 
\end{Claim}

{\em Proof}.  First we will prove that $Y(\mathbf{a}, 0)$ has only isolated singularity at $\bar{0}$. 
For a point $p \in \mu^{-1}(0)$, let us denote by $Op$ the $\mathbf{C}^*$-orbit 
containing $p$. It is easily checked that the closure $\overline{Op}$ of $Op$ contains the 
origin $0$ if and only if $z_i(p) = 0$ for all $i$ or $w_i(p) = 0$ for all $i$. 
Define $$F := \{z_1 = ... = z_n = 0\} \:\: G := \{w_1 = ... = w_n = 0\}.$$ 
Both $F$ and $G$ are $n$-dimensional affine spaces contained in $\mu^{-1}(0)$. 
Assume that $p \notin F \cup G$. Let us consider the stabilizer group $\mathbf{C}^*_p \subset 
\mathbf{C}^*$ of $p$. 
Since $p \notin F$, $z_i(p) \ne 0$ for some $i$. 
Similarly, $w_j(p) \ne 0$ for some $j$ because $p \notin G$. If one can choose $i$ and $j$ so that 
$i \ne j$, then $t \in \mathbf{C}^*_p$ must satisfy $t^{a_i} = 1$ and $t^{-a_j} = 1$, which means $t = 1$ because $a_i$ and $a_j$ are 
coprime. 
Suppose, to the contrary, that there is a unique $i$ such that $z_i(p) \ne 0$ and $w_i(p) \ne 0$ and 
$z_j(p) = w_j(p) = 0$ for any other $j$. Then 
$$0 = \mu(p) = a_iz_i(p)w_i(p), $$
which means that $z_i(p) = 0$ or $w_i(p) = 0$. But, this is absurd. Hence, such a case does not occur. 
As a consequence, we have $\mathbf{C}^*_p = 1$ for $p \in \mu^{-1}(0) - F - G$.     
Here let us consider the quotient map 
$$\pi: \mu^{-1}(0) \to Y(\mathbf{a}, 0).$$
Then $\pi^{-1}(\bar{0}) = F \cup G$ and $$\pi\vert_{\mu^{-1}(0) - \pi^{-1}(\bar{0})}: \mu^{-1}(0) - \pi^{-1}(\bar{0}) \to 
Y(\mathbf{a}, 0) - \{\bar{0}\}$$ is a $\mathbf{C}^*$-bundle.  
Since $\mu^{-1}(0)$ is smooth outside the origin $0$, $\mu^{-1}(0) - \pi^{-1}(\bar{0})$ is also smooth. This means 
that $Y(\mathbf{a}, 0) - \{\bar{0}\}$ is smooth. 
As we will see below Claim \ref{2}, there is a partial resolution $\nu^+: Y(\mathbf{a}, +) \to Y(\mathbf{a}, 0)$ with 
only quotient singularities. The symplectic form $\omega$ on $\mathbf{C}^{2n}$ descends to a symplectic form $\omega^+$ (resp. 
$\omega^0$) on $Y(\mathbf{a}, +)_{reg}$ (resp. $Y(\mathbf{a}, 0)_{reg}$). Moreover, $\omega^+ = (\nu^+)^*\omega^0$. 
The partial resolution $Y(\alpha, +)$ has symplectic quotient singularities by \cite{Be}. Therefore $Y(\alpha, 0)$ also has symplectic singularities.   
$\square$

\begin{Claim}\label{3} We have $\pi_1(Y(\mathbf{a}, 0) - \{\bar{0}\}) = 1$. 
\end{Claim}    

{\em Proof}. By the proof of Claim \ref{2}, $\pi\vert_{\mu^{-1}(0) - \pi^{-1}(\bar{0})}$ is a $\mathbf{C}^*$-bundle. Then there is a 
surjection $\pi_1(\mu^{-1}(0) - \pi^{-1}(\bar{0})) \to 
\pi_1(Y(\mathbf{a}, 0) - \{\bar{0}\})$ by the homotopy exact sequence. Therefore it suffices to show that $\pi_1(\mu^{-1}(0) - \pi^{-1}(\bar{0})) = 1$. 
The central fiber $\mu^{-1}(0)$ is a hypersurface of $\mathbf{C}^{2n}$ with an isolated singularity at $0$. 
Note that $\dim \mu^{-1}(0) = 2n-1$. By [\cite{M}, Theorem 5.2] the local fundamental group of 
an isolated hypersurface singularity of dim $\geq 3$ is trivial.
Moreover,  
$\mu^{-1}(0)$ is stable under the scaling $\mathbf{C}^*$-action 
$$(z_1, ..., z_n, w_1, ..., w_n) \to (tz_1, ..., tz_n, tw_1, ..., tw_n).$$           
Then, by [ibid] together with this $\mathbf{C}^*$-action, we have 
$\pi_1(\mu^{-1}(0) - \{0\}) = 1$. 
On the other hand, 
$$\mathrm{Codim}_{\mu^{-1}(0)} F \cup G = n - 1 \geq 2.$$ 
Since $\mu^{-1}(0) - \pi^{-1}(\bar{0})$ is obtained from the complex manifold $\mu^{-1}(0) - \{0\}$ 
by removing a closed subspace at least codimension 2, we see that 
$$\pi_1(\mu^{-1}(0) - \pi^{-1}(\bar{0})) \cong \pi_1(\mu^{-1}(0) - \{0\}) = 1.$$ 
$\square$ 
     
For a character $\alpha \in \mathrm{Hom}_{alg.gp}(\mathbf{C}^*, \mathbf{C}^*)$, a polynomial function $f$ on $\mathbf{C}^{2n}$ 
is called $\alpha$-invariant if $f(tp) = \alpha(t)f(p)$ holds for $t \in \mathbf{C}^*$ and $p \in 
\mathbf{C}^{2n}$. A point $p \in \mathbf{C}^{2n}$ is called $\alpha$-semsistable if there is a positive integer $l > 0$ 
and an $l\cdot \alpha$-invariant polynomial function $f$ such that $f(p) \ne 0$.  
We put $$\mu^{-1}(0)^{\alpha-ss} := \{p \in \mu^{-1}(0)\vert \:\: p \: \mathrm{is}\: \alpha \mathrm{-semistable}\}.$$ 
and take its GIT-quotient $Y(\mathbf{a}, \alpha)$. By definition, there is a quotient map 
$$\pi^{\alpha}: \mu^{-1}(0)^{\alpha-ss} \to Y(\mathbf{a}, \alpha).$$  
One can write $\alpha$ as $\alpha(t) = t^{m(\alpha)}$ with a unique integer $m(\alpha)$. 
$\mu^{-1}(0)^{\alpha-ss}$ changes according as $m(\alpha) > 0$, $m(\alpha) = 0$ or $m(\alpha) < 0$. 
Hence we write $\mu^{-1}(0)^+$ for $\mu^{-1}(0)^{\alpha-ss}$ when $m(\alpha) > 0$, and 
write $\mu^{-1}(0)^-$ for $\mu^{-1}(0)^{\alpha-ss}$ when $m(\alpha) < 0$. Similarly we respectively denote by  
$Y(\mathbf{a}, +)$ and $Y(\mathbf{a}, -)$ their quotients. 

Now one can check 
$$ \mu^{-1}(0)^+ = \mu^{-1}(0) - F, \:\:\: \mu^{-1}(0)^- = \mu^{-1}(0) - G.$$
The inclusion maps $\mu^{-1}(0)^+ \subset \mu^{-1}(0)$ and 
$\mu^{-1}(0)^- \subset \mu^{-1}(0)$ respectively induce maps $\nu^+: Y(\mathbf{a}, +) \to Y(\mathbf{a}, 0)$ and 
$\nu^-: Y(\mathbf{a}, -) \to Y(\mathbf{a}, 0)$. 
It is easily checked that $\nu^+$ and $\nu^-$ are both isomorphisms over $Y(\mathbf{a}, 0) - \{\bar{0}\}$ and    
$(\nu^+)^{-1}(\bar{0}) \cong \mathbf{P}(a_1, ..., a_n)$, $(\nu^-)^{-1}(\bar{0}) \cong \mathbf{P}(a_1, ..., a_n)$.
Both $Y(\mathbf{a}, +)$ and $Y(\mathbf{a}, -)$ are crepant partial resolutions of $Y(\mathbf{a}, 0)$ with only isolated quotient singularities. 
Therefore, $\nu^+$ and $\nu^{-}$ are (projective) $\mathbf{Q}$-factorial terminalizations 
of $Y(\mathbf{a},0)$.
The diagram 
$$ Y(\mathbf{a}, +) \stackrel{\nu^+}\rightarrow Y(\mathbf{a}, 0) \stackrel{\nu^-}\leftarrow Y(\mathbf{a}, -)$$ 
is a flop in the sense that the proper transform $L'$ of a $\nu^+$-negative line bundle $L$ is $\nu^{-}$-ample. 
Since the relative Picard numbers $\rho(\nu^+)$ and $\rho(\nu^{-})$ are both equal to $1$, we see that 
$Y(\mathbf{a}, 0)$ has exactly two different {\bf Q}-factorial terminalizations. Assume that $\mathbf{a'}$ is an $n$-tuple 
of poisitive integers different from 
$\mathbf{a}$. Then $Y(\mathbf{a}, 0) \ne Y(\mathbf{a}', 0)$. In fact, if they coincide, then their $\mathbf{Q}$-factorial terminalizations 
must coincide.  In particular, their exceptional loci coincide; but $\mathbf{P}(a_1, ..., a_n) \ne \mathbf{P}(a'_1, ..., a'_n)$, which is a contradiction.  As a consequence, we have Proposition \ref{1}. 
\vspace{0.2cm}

{\bf Maximal weight}. Let $(X, \omega)$ be a conical symplectic variety with coordinate ring $R$. By definition, 
$R$ is positively graded $R = \oplus_{i \ge 0}R_i$ with $R_0 = \mathbf{C}$ and $\omega$ is a homogeneous symplectic form on 
$X_{reg}$. The $\mathbf{C}^*$-action on $X$ determined by the grading is called a conical $\mathbf{C}^*$-action.  
Let $x_1, ..., x_k \in R$ be minimal homogeneous generators of $R$. Put $l_i := \mathrm{deg}(x_i)$. We assume that the conical 
$\mathbf{C}^*$-action is effective. In other words, $\mathrm{GCD}(l_1, ..., l_k) = 1$.   
We define $$l := \mathrm{max}\{l_1, ..., l_k\}$$ and call it the {\em maximal weight} of $X$. 
If we fix the dimension and the maximal weight, then there are only finitely many conical symplectic varieties with such properties (\cite{Na 1}). 
For example, if a nilpotent orbit $O$ of a complex semisimple Lie algebra $\mathfrak{g}$ has normal closure $\bar{O}$, then 
$\bar{O}$ is a conical symplectic variety with maximal weight $1$. Conversely, any conical symplectic variety $X$ with maximal 
weight $1$, has this form (\cite{Na 2}). Note that, when $\bar{O}$ is not normal, its normalization $\tilde{O}$ may possibly have 
maximal weight $> 1$.  

We shall explain here that our toric hyperk\"{a}hler construction gives an example with an arbitrary maximal weight. 

Assume that the $n$-tuple $\mathbf{a}$ is of the form $(1, ..., 1, 2l -1)$ with a positive integer $l$. It is easily checked that the $\mathbf{C}$-algebra  
$R := \mathbf{C}[Y(\mathbf{a}, 0)]$ is generated by the monomials\footnote{There is a relation: $z_1w_1 + \cdots + z_{n-1}w_{n-1} + mz_nw_n 
= 0$.} 
$$z_iw_j \:\: (1 \le i \le n-1, \: 1 \le j \le n-1), \:\: z_nw_n$$
$$z_1^{d_1}\cdots z_{n-1}^{d_{n-1}}w_n, \:\: z_nw_1^{d_1}\cdots w_{n-1}^{d_{n-1}} \:\: (d_1 + ... + d_{n-1} = 2l -1).$$
For the scaling $\mathbf{C}^*$-action $(z_1, ..., z_n, w_1, ..., w_n) \to (tz_1, ..., tz_n, tw_1, ..., tw_n)$, 
the monomials on the first row have weight $2$ and those on the second row have weight $2l$.
As it stands, the scaling $\mathbf{C}^*$-action is not effective. Now one can introduce a new $\mathbf{C}^*$-action on $\mathbf{C}[Y(\mathbf{a}, 0)]$ so that the monomials on the first row have weight $1$ and those on the second row have weight $l$. 
Then, together with the new $\mathbf{C}^*$-action, $(Y(\mathbf{a}, 0), \omega)$ is a conical symplectic variety of dimension $2n-2$ with maximal weight $l$. \vspace{0.2cm}

\section{Further generalization}

The examples above can be generalized as follows. Let $n$ and $d$ be positive integers 
such that $n \geq d +2$. 
Let $A$ (resp. $B$) be a $d \times n$-matrix (resp. an $n \times (n-d)$-matrix) with integer coefficients so that  
the following sequence is exact:   
$$0 \to \mathbf{Z}^{n-d} \stackrel{B}\to \mathbf{Z}^n \stackrel{A}\to \mathbf{Z}^d \to 0.$$
Assume the following two conditions: 

(*) Any $d \times d$-minor matrix of $A$ has non-zero determinant.  
\vspace{0.2cm}
 
(**) Take an arbitrary  $d \times (d+1)$-minor matrix $A'$ of $A$,  and consider all $d \times d$-minor matrices contained in 
$A'$. Then the greatest common divisor of their determinants is $1$.      
\vspace{0.2cm}

When $d = 1$, they are nothing but the conditions imposed in the previous example.  
A $d$ dimensional algebraic torus $(\mathbf{C}^*)^d$ acts on $\mathbf{C}^{2n}$ by 
$$ (z_1, ..., z_n, w_1, ..., w_n) \to$$ $$(t_1^{a_{11}}t_2^{a_{21}}\cdot \cdot \cdot t_d^{a_{d1}}z_1, ..., t_1^{a_{1n}}t_2^{a_{2n}}\cdot \cdot \cdot t_d^{a_{dn}}z_n, t_1^{-a_{11}}t_2^{-a_{21}}\cdot \cdot \cdot t_d^{-a_{d1}}w_1, ...,    
t_1^{-a_{1n}}t_2^{-a_{2n}}\cdot \cdot \cdot t_d^{-a_{dn}}w_n).$$
Define a symplectic form  $\omega$ on $\mathbf{C}^{2n}$ by $$\omega := \sum_{1 \le i \le n} dw_i \wedge dz_i.$$  
Then the $(\mathbf{C}^*)^d$-action on $(\mathbf{C}^{2n}, \omega)$ is a Hamiltonian action and 
its moment map $\mu$ is given by 
$$\mathbf{C}^{2n} \to \mathbf{C}^d\:\:\: (z_1, ..., z_n, w_1, ..., w_n) \to (\sum_{1 \le j \le n} a_{1j}z_jw_j, ..., \sum_{1 \le j \le n} a_{dj}z_jw_j). 
$$
Define $$Y(A, 0) := \mu^{-1}(\mathbf{0})/\hspace{-0.1cm}/(\mathbf{C}^*)^d$$ and call it a toric hyperk\"{a}hler variety associated with $A$. 
Let $$\pi: \mu^{-1}(\mathbf{0}) \to Y(A, 0)$$ be the quotient map. We put $\bar{0} := \pi(0)$ for $0 \in \mathbf{C}^{2n}$.  

\begin{Prop}\label{4}  The toric hyperk\"{a}hler variety $Y(A, 0)$ is a symplectic variety of dimension $2(n-d)$ with an isolated singularity and $\pi_1(Y(A, 0)_{reg}) = 1$.
\end{Prop}

We will prove the first part of Proposition \ref{4} in Lemma \ref{5}, and the second part in Lemma \ref{10}.  

\begin{Lem}\label{5} The variety $Y(A,0)$ is a symplectic variety of dimension $2(n-d)$ with an isolated singularity at $\bar{0}$.
\end{Lem}

{\em Proof}. Lemma \ref{5} will be proved after Claim \ref{6}, Claim \ref{7}, Claim \ref{8} and Corollary \ref{9}. Put $f_i := \sum a_{ij}z_jw_j$ for $1 \le i \le d$. Then the Jacobi matrix  
$$ J(\mu) :=  \left( \begin{array}{cccccc} 
\frac{\partial f_1}{\partial z_1} & \cdots & \frac{\partial f_1}{\partial z_n} & \frac{\partial f_1}{\partial w_1} & \cdots & \frac{\partial f_1}{\partial w_n}\\ 
\cdots & \cdots & \cdots & \cdots & \cdots & \cdots \\ 
\frac{\partial f_d}{\partial z_1} & \cdots & \frac{\partial f_d}{\partial z_n} & \frac{\partial f_d}{\partial w_1} & \cdots & \frac{\partial f_d}{\partial w_n} 
\end{array} \right)  $$ is given by 
$$\left( \begin{array}{cccccc} 
a_{11}w_1 & \cdots & a_{1n}w_n & a_{11}z_1 & \cdots & a_{1n}z_n\\ 
\cdots & \cdots & \cdots & \cdots & \cdots & \cdots \\ 
a_{d1}w_1 & \cdots & a_{dn}w_n & a_{d1}z_1 & \cdots & a_{dn}z_n
\end{array} \right).$$ 
$\mathrm{Sing}(\mu^{-1}(\mathbf{0}))$ consists of the points of $\mu^{-1}(\mathbf{0})$ where 
$d \times d$-minors of $J(\mu)$ are all zero.  For a subset $I := \{i_1, ..., i_d\}$ of $\{1, ..., n\}$,  
we put    
$$A_I := \left( \begin{array}{ccc} 
a_{1,i_1} & \cdots & a_{1, i_d}\\ 
\cdots & \cdots & \cdots \\ 
a_{d, i_1} & \cdots & a_{d, i_d}
\end{array} \right).$$ 
Then a $d \times d$ minor of $J(\mu)$ is of the form 
$$\pm \mathrm{det}(A_I) z^Jw^K$$ 
for a subset $I$ of $\{1, ..., n\}$ with $\vert I \vert = d$ and its decomposition $I = J \sqcup K$.  
Since $\mathrm{det}(A_I) \ne 0$ by condition (*), we see that 
$$\mathrm{Sing}(\mu^{-1}(\mathbf{0})) = \bigcap_{I \subset \{1, ..., n\}, \: \vert I \vert =d, \: I = J \sqcup K}\{z^Jw^K = 0\} 
\cap \mu^{-1}(\mathbf{0}).$$
Let us look more closely at a point $$p \in \bigcap_{I \subset \{1, ..., n\}, \: \vert I \vert =d, \: I = J \sqcup K}\{z^Jw^K = 0\}.$$
Define  
$$I(p) := \{i \: \vert \: z_i(p) \ne 0 \: \mathrm{or} \: w_i(p) \ne 0\}.$$
If $\vert I(p) \vert \geq d$, then take a subset $I$ of $I(p)$ such that $\vert I \vert = d$. 
By the definition of $I(p)$ one can find a decomposition $I = J  \sqcup K$ such that $z^J(p)w^K(p) \ne 0$, which is a contradiction. 
Hence, one can choose a sequence of integers $1 \le i_1 < \cdots < i_{d-1} \le n$ such that 
$$z_i(p) = w_i(p) = 0 \:\: \mathrm{for}\:\:  i \notin \{i_1, ..., i_{d-1}\}.$$
Assume in addition that $p \in  \mu^{-1}(\mathbf{0})$. Then 
$$\mathbf{a}_{i_1}z_{i_1}(p)w_{i_1}(p) + \cdots + \mathbf{a}_{i_{d-1}}z_{i_{d-1}}(p)w_{i_{d-1}}(p) = \mathbf{0},$$
where $\mathbf{a}_i$ denotes the $i$-th column of $A$. By (*), the column vectors 
$\mathbf{a}_{i_1}$, ..., $\mathbf{a}_{i_{d-1}}$ are linearly independent; hence we have 
$$z_{i_1}(p)w_{i_1}(p) = \cdots =  z_{i_{d-1}}(p)w_{i_{d-1}}(p) = 0.$$

{\bf Remark}. By the observation above, we see that $\dim \mathrm{Sing}(\mu^{-1}(\mathbf{0})) = d-1$.

\begin{Claim}\label{6} If $p \in \mathrm{Sing}(\mu^{-1}(\bar{0}))$, then $p \in \pi^{-1}(\bar{0})$. 
\end{Claim}

{\em Proof}. Let $$\lambda: \mathbf{C}^* \to (\mathbf{C}^*)^d, \:\: t \to t_1^{\lambda_1}t_2^{\lambda_2}\cdots t_d^{\lambda_d} \: (\lambda_1, ..., \lambda_d \in \mathbf{Z})$$ be a 1-parameter subgroup of $(\mathbf{C}^*)^d$. 
Then $\mathbf{C}^*$ acts on $\mu^{-1}(\mathbf{0})$ via $\lambda$. 
It suffices to show that $\lim_{t \to 0}\lambda (t)\cdot p = \mathbf{0}$ for a suitable $\lambda$. 
For simplicity, we only consider the case where $(i_1, ..., i_{d-1}) = (1, ..., d-1)$ and $w_1(p) = \cdots = w_{d-1}(p) = 0$ because 
other cases can be treated in the same way. In this case 
$$\lambda (t)\cdot p = (t^{a_{1,1}\lambda_1 + \cdots + a_{d,1}\lambda_d}z_1(p), ..., 
t^{a_{1, d-1}\lambda_1 + \cdots + a_{d, d-1}\lambda_d}z_{d-1}(p), 0, ..., 0).$$  
Since $\mathbf{a}_1$, ..., $\mathbf{a}_{d-1}$ are linearly independent vectors of $\mathbf{R}^d$, one can choose 
$(\lambda_1, ..., \lambda_{d-1})$ so that  
$$a_{1,j}\lambda_1 + \cdots + a_{d,j}\lambda_d > 0, \:\: (j = 1, ..., d-1).$$ 
Then $\lim_{t \to 0}\lambda (t)\cdot p = \mathbf{0}$ for such $\lambda$.  $\square$   
\vspace{0.2cm}

For the $(\mathbf{C}^*)^d$-action on $\mu^{-1}(\mathbf{0})$, denote by $(\mathbf{C}^*)^d_p$ the stabilizer group of $p$.

\begin{Claim}\label{7} For a point $p \in \mu^{-1}(\mathbf{0})$, the following conditions are equivalent: 

(1) $p \in \mu^{-1}(\mathbf{0})_{reg}$.

(2) {\em $(\mathbf{C}^*)^d_p$ is a finite group.} 
\end{Claim}

{\em Proof}. The $(\mathbf{C}^*)^d$-action on $\mathbf{C}^{2n}$ induces a vector $\zeta_a \in T_p\mathbf{C}^{2n}$ 
for each $a \in \mathrm{Lie}((\mathbf{C}^*)^d)$. We have $a \in \mathrm{Lie}((\mathbf{C}^*)^d_p)$ if and only if 
$\zeta_a = 0$. Note that $\zeta_a = 0$ if and only if $\omega_p(v, \zeta_a) = 0$ for all $v \in T_p\mathbf{C}^{2n}$. 
Since $\mu$ is the moment map, we have 
$$\omega_p(v, \zeta_a) = \langle d\mu_p(v), a \rangle.$$ Hence, 
$\zeta_a = 0$ if and only if $\langle d\mu_p(v), a \rangle = 0$ for all $v \in T_p(\mathbf{C}^{2n})$. 
This means that $d\mu_p$ is surjective if and only if $\mathrm{Lie}((\mathbf{C}^*)^d_p) = 0$. 
This is nothing but the claim. $\square$


\vspace{0.2cm}

\begin{Claim}\label{8} Assume that $(\mathbf{C}^*)_p$ is a nontrivial finite group for $p \in \mu^{-1}(\mathbf{0})$. Then 
$p \in \pi^{-1}(\bar{0})$. 
\end{Claim}

{\em Proof}. 
Put $$I(p) := \{i \: \vert \: z_i(p) \ne 0 \: \mathrm{or} \: w_i(p) \ne 0\}.$$ 
We first show that  $\vert I(p) \vert = d$. In fact, if $\vert I(p) \vert < d$, then $(\mathbf{C}^*)^d_p$ must be infinite. 
If $\vert I(p) \vert \geq d+1$, then $(\mathbf{C}^*)^d_p = 1$ by condition (**). Therefore, we have $\vert I(p) \vert = d$. 
Write $I(p) = \{i_1, ..., i_d\}$. Since $p \in \mu^{-1}(\mathbf{0})$, we have
$$\mathbf{a}_{i_1}z_{i_1}(p)w_{i_1}(p) + \cdots + \mathbf{a}_{i_{d}}z_{i_{d}}(p)w_{i_{d}}(p) = \mathbf{0}.$$  
By condition (*), the column vectors $\mathbf{a}_{i_1}$, ..., $\mathbf{a}_{i_d}$ are basis of $\mathbf{R}^d$; hence we have  
$$z_{i_1}(p)w_{i_1}(p) = \cdots =  z_{i_{d}}(p)w_{i_{d}}(p) = 0.$$
The following argument is quite similar to Claim \ref{6}. As in Claim \ref{6}, we have to prove that $\lim_{t \to 0}\lambda (t)\cdot p = \mathbf{0}$ 
for a suitable 1-parameter subgroup $\lambda: \mathbf{C}^* \to (\mathbf{C}^*)^d$. 
For simplicity, we only consider a typical case  
where $(i_1, ..., i_d) = (1, ..., d)$ and $w_1(p) = \cdots w_d(p) = 0$. 
In this case 
$$\lambda (t)\cdot p = (t^{a_{1,1}\lambda_1 + \cdots + a_{d,1}\lambda_d}z_1(p), ..., 
t^{a_{1, d}\lambda_1 + \cdots + a_{d, d}\lambda_d}z_d(p), 0, ..., 0).$$  
Since $\mathbf{a}_1$, ..., $\mathbf{a}_d$ are basis of $\mathbf{R}^d$, one can choose 
$(\lambda_1, ..., \lambda_d)$ so that  
$$a_{1,j}\lambda_1 + \cdots + a_{d,j}\lambda_d > 0, \:\: (j = 1, ..., d).$$ 
Then $\lim_{t \to 0}\lambda (t)\cdot p = \mathbf{0}$ for such $\lambda$. $\square$ 
\vspace{0.2cm}

\begin{Cor}\label{9} If $p \in \mu^{-1}(\mathbf{0}) - \pi^{-1}(\bar{0})$, then $(\mathbf{C}^*)^d_p = 1$.  
\end{Cor}

{\em Proof}. By Claim \ref{6}, $p \in \mu^{-1}(\mathbf{0})_{reg}$. Then, by Claim \ref{7}, $(\mathbf{C}^*)^d_p$ is a finite group. 
If $(\mathbf{C}^*)^d_p$ is nontrivial, then $p \in \pi^{-1}(\bar{0})$. Therefore, $(\mathbf{C}^*)^d_p = 1$. 
$\square$ 
\vspace{0.2cm}

We are now in a position to prove Lemma \ref{5}.  Let us consider the commutative diagram 
\begin{equation} 
\begin{CD}
\mu^{-1}(\mathbf{0}) - \pi^{-1}(\bar{0}) @>>> \mu^{-1}(\mathbf{0}) \\
@V{\pi\vert_{\pi^{-1}(Y(A, 0) - \{\bar{0}\})}}VV  @V{\pi}VV \\ 
Y(A, 0) - \{\bar{0}\}  @>>> Y(A, 0).        
\end{CD} 
\end{equation}
Here $\mu^{-1}(\mathbf{0}) - \pi^{-1}(\bar{0})$ is nonsingular by Claim \ref{6}. 
By Corollary \ref{9},  
 $\pi\vert_{\pi^{-1}(Y(A, 0) - \{\bar{0}\})}$ is a $(\mathbf{C}^*)^d$-bundle. 
This means that $Y(A, 0) - \{\bar{0}\}$ is nonsingular.  Since $\dim \mu^{-1}(\mathbf{0}) = 2n -d$, we have 
$\dim Y(A, 0) = 2(n-d)$.  We finally prove that $Y(A, 0)$ has symplectic singularities. This is already proved in 
\cite{B-F}, Proposition 5.1 and \cite{Bel}, Proposition 2.5. But we include here the proof for the sake of completeness.  
For a character $\alpha \in \mathrm{Hom}_{alg.gp}((\mathbf{C}^*)^d, \mathbf{C}^*)$, we define 
$Y(A, \alpha)$ as the GIT quotient of $\mu^{-1}(0)^{\alpha-ss}$ by $(\mathbf{C}^*)^d$. 
By the proof of [\cite{HS}, Proposition 6.2], $Y(A, \alpha)$ has only quotient singularities if $\alpha$ is generic. Moreover, 
the symplectic 2-form $\omega$ on $\mathbf{C}^{2n}$ descends to a symplectic form $\omega^{\alpha}$ (resp. 
$\omega^0$) on $Y(A, \alpha)_{reg}$ (resp. $Y(A, 0)_{reg}$).  
There is a birational projective morphism $\nu^{\alpha}: Y(A, \alpha) \to Y(A, 0)$ and $(\nu^{\alpha})^*\omega^0 
= \omega^{\alpha}$. Since $Y(A, \alpha)$ has symplectic quotient singularities by \cite{Be}, we see that 
$Y(A, 0)$ also has symplectic singularities.  
$\square$ 
\vspace{0.2cm}

\begin{Lem}\label{10} When $n \geq d + 2$, we have $\pi_1(Y(A, 0)_{reg}) = 1$. 
\end{Lem}

{\em Proof}. Since $\dim \mathrm{Sing}(\mu^{-1}(\mathbf{0})) = d-1$ by Remark above Claim \ref{6}, we have 
$$\mathrm{Codim}_{\mu^{-1}(\mathbf{0})}\mathrm{Sing}(\mu^{-1}(\mathbf{0})) = (2n-d) - (d-1) = 2(n-d) + 1 \geq 5.$$
In particular, $\mu^{-1}(\mathbf{0})$ has complete intersection singularities with $\mathrm{Codim}_{\mu^{-1}(\mathbf{0})}\mathrm{Sing}(\mu^{-1}(\mathbf{0}))
\geq 3$. Then we conclude that $\pi_1(\mu^{-1}(\mathbf{0})_{reg}) = 1$ by [\cite{HL}, Corollary 3.2.2] in the following way.  First recall the notion of 
{\em rectified homotopical depth}. Let $X$ be a 
reduced complex analytic space and take a Whitney stratification of $X$. Let $X_i$ be the union of the strata of 
dimension $\le i$. The rectified homotopical depth of $X$ is $\geq m$ if, for any point $x \in X_i - X_{i-1}$, there is a 
sufficiently small contractible open neghborhood $U$ of $x$ in $X$ such that 
$\pi_k(U, U - X_i) = 1$ for $k < m -i$. By [ibid], the rectified homotopical depth of $X$ is equal to $\dim X$ if 
$X$ is locally of complete intersection. We apply this result to $X: = \mu^{-1}(\mathbf{0})$ and prove that 
$\pi_1(\mu^{-1}(\mathbf{0})_{reg}) \cong \pi_1(\mu^{-1}(\mathbf{0}))$. Note that $X_{d-1} = \mathrm{Sing}(\mu^{-1}(\mathbf{0}))$.
It suffices to prove that a connected covering $Z_{d-1} \to X - X_{d-1}$ 
uniquely extends to a connected covering $Z \to X$.  We prove that a connected covering $Z_{d-i} \to X - X_{d-i}$ uniquely 
extends to a connected covering $Z_{d-i-1} \to X - X_{d-i-1}$ by the induction on $i$. 
For  $x \in X_{d-i} - X_{d-i-1}$, we take a contractible open neighborhood $U$ of $x \in X$. Since the rectified homotopical depth 
of $X$ is $2n - d$, we have $\pi_k(U, U - X_{d-i}) = 1$ for $k < 2n-d - (d - i ) = 2(n-d) + i  \le 4 + i$. 
In particular, $\pi_1(U, U - X_{d-i}) = \pi_2(U, U - X_{d-i}) 
= 1$.   By the homotopy exact sequence
$$\pi_2(U, U - X_{d-i}) \to \pi_1(U - X_{d-i}) \to \pi_1(U) \to \pi_1(U, U- X_{d-i})$$ 
we have an isomorphism $\pi_1(U - X_{d-i}) \cong \pi_1(U)$. This means that $Z_{d-i} \to X - X_{d-i}$ uniquely 
extends to a connected covering $Z_{d-i-1} \to X - X_{d-i-1}$. Finally we get a connected covering 
$Z \to X$.  We next prove that  $\pi_1(\mu^{-1}(\mathbf{0})) = 1$. In fact, $\mu^{-1}(\mathbf{0})$ has 
a scaling $\mathbf{C}^*$ action defined by 
$$(z_1, ..., z_n, w_1, ..., w_n) \to (tz_1, ..., tz_n, tw_1, ..., tw_n), \:\:\: t \in \mathbf{C}^*.$$ 
By using this $\mathbf{C}^*$-action, $\mu^{-1}(\mathbf{0})$ is homotopy equivalent to a sufficiently small contractible 
open neighborhood of $\mathbf{0} \in \mu^{-1}(\mathbf{0})$; hence, $\pi_1(\mu^{-1}(\mathbf{0)}) = 1$.  As a consequence, 
we have $\pi_1(\mu^{-1}(\mathbf{0})_{reg}) = 1$. 

We next prove that $\pi_1(\mu^{-1}(\mathbf{0}) - \pi^{-1}(\bar{0})) = 1$. By Claim \ref{6}, $\mu^{-1}(\mathbf{0}) - \pi^{-1}(\bar{0}) = 
\mu^{-1}(\mathbf{0})_{reg} - \pi^{-1}(\bar{0})$. Note that 
$$\pi^{-1}(\bar{0}) \subset \{z_iw_i = 0\: \vert \: 1 \le \forall i \le n\}$$ 
and the right hand side has dimension $n$. Hence $\dim \pi^{-1}(\bar{0}) \leq n$ and 
$$\mathrm{Codim}_{\mu^{-1}(0)}\pi^{-1}(\bar{0}) \geq (2n-d) - n = n - d \geq 2.$$ 
Therefore, we have an  isomorphism
$$\pi_1(\mu^{-1}(\mathbf{0})_{reg} - \pi^{-1}(\bar{0})) \cong \pi_1(\mu^{-1}(\mathbf{0})_{reg}) ( = 1)$$

Since $$\pi\vert_{\pi^{-1}(Y(A, 0) - \{\bar{0}\})}: \mu^{-1}(\mathbf{0}) - \pi^{-1}(\bar{0}) \to Y(A, 0) - \{\bar{0}\}$$ 
is a $(\mathbf{C}^*)^d$-bundle, there is a surjection 
$$\pi_*: \pi_1(\mu^{-1}(\mathbf{0}) - \pi^{-1}(\bar{0})) \to \pi_1(Y(A, 0) - \{\bar{0}\}).$$
This means that $\pi_1(Y(A, 0) - \{\bar{0}\}) = 1$. $\square$ 

\section{Contact structures}

Our examples $Y(A, 0)$ are all non-{\bf Q}-factorial because they have non-trivial projective {\bf Q}-factorial terminalizations.
Towards the classification of symplectic singularities, it would be interesting to find new examples of  
{\bf Q}-factorial conical symplectic varieties with an isolated singularity whose local fundamental group is trivial.   

As explained below, this problem is also closely related with the LeBrun-Salamon conjecture for a 
contact Fano manifold. 

A complex manifold $M$ of dimension $2n-1$ is called a contact complex manifold if 
there exists a holomorphic line bundle $L$ on $M$ and a surjection $\theta: \Theta_M \to L$ such that the 
following condition hold: 

If we define $F := \mathrm{Ker}(\theta)$, then 
$$F \times F \to L, \:\: (x,y) \to \theta([x,y])$$ is a nondegenerate skew symmetric form. 

Here $x$ and $y$ are regarded as sections of $\Theta_M$ and $[x, y]$ is the bracket product of vector fields. 
Such a line bundle $L$ is called a contact line bundle. Then the canonical line bundle $K_M$ of $M$ satisfies 
$K_M \cong L^{-n}$. In particular, if $M$ is a Fano manifold, then $L$ is an ample line bundle. 
The LeBrun-Salamon conjecture asserts that, if $M$ is a contact Fano manifold, then 
$M$ is isomorphic to the projectivization $\mathbf{P}(O_{min})$ of the minimal nilpotent orbit $O_{min}$ of a complex simple Lie algebra 
$\mathfrak{g}$. The conjecture is true if  $\dim M \leq 5$ or $b_2(M) > 1$ (cf. \cite{Ye}, \cite{Dr}, \cite{KPSW}).  

Let $(M, L)$ be a contact Fano manifold. Then the contact structure on $M$ induces a symplectic form 
$\omega_{(L^{-1})^{\times}}$ on $(L^{-1})^{\times} := L^{-1} - (0\mathrm{-section})$. Since $L^{-1}$ is negative, we can contract the $0$-section of $L^{-1}$ to get an affine variety $X$ with an isolated singularity $0$. Then $\omega_{(L^{-1})^{\times}}$ determines 
a symplectic form $\omega$ on $X - \{0\}$.  Since $L^{-1}$ has a natural $\mathbf{C}^*$-action, 
$X$ also has a $\mathbf{C}^*$-action and $(X, \omega)$ becomes a conical symplectic variety with $wt(\omega) = 1$. 
By construction, $\mathbf{C}^*$ acts on $X - \{0\}$ freely.  
Conversely, assume that  $X$ is a conical symplectic variety with an isolated 
singularity at $0 \in X$ such that the $\mathbf{C}^*$-action is {\em free} on $X - \{0\}$. Embed $X$ in an affine space $\mathbf{C}^k$ 
by using the minimal homogeneous generators of the coordinate ring $R$ of $X$. Denote by  $l_1, ..., l_k$ the weights of the 
generators. Let $f: Z \to X$ be the weighted blowing-up of $X$ with respect to $(l_1, ..., l_k)$. Then $M := f^{-1}(0)$ is 
a contact Fano manifold (cf. [\cite{Na 1}, \S 2], [\cite{Na 3}, Lemma 2.1]). As the LeBrun-Salamon conjecture holds true if $b_2(M) > 1$ or $\dim M \leq 5$, we are interested in 
the case $b_2(M) = 1$ and $\dim M \geq 7$. Then, since $f$ is a resolution such that $\mathrm{Exc}(f)$ has 2-nd Betti number 
$1$, $X$ must be {\bf Q}-factorial with 
$\dim X \geq 8$. It would be a challenging problem to ask if there exists a new {\bf Q}-factorial conical symplectic variety $X$, which 
satisfies 

(i)  $X$ is not $\bar{O}_{min}$, and $\dim X \geq 8$, 

(ii) the origin $0 \in X$ is an isolated singularity, 

(iii) $\pi_1(X - \{0\}) = 1$, and finally 

(iv) the conical $\mathbf{C}^*$-action on $X$ is free outside $0 \in X$.

\begin{center}
Research Institute for Mathematical Sciences, Kyoto University, Oiwake-cho, Kyoto, Japan

E-mail address: namikawa@kurims.kyoto-u.ac.jp  
\end{center}


\begin{thebibliography}{99} 


\bibitem[Be]{Be} Beauville, A.: Symplectic singularities, Invent. Math. {\bf 139} (2000), no.3, 541-549


\bibitem[Bel]{Bel} Bellamy, G.: Coulomb branches have symplectic singularities. Lett. Math. Phys. {\bf 113} (2023), no. 5, Paper No. 104, 8 pp.

\bibitem[BBFJLS]{BBFJLS} Bellamy, G., Bonnaf\'{e}, C.,  Fu, B., Juteau, D.,  Levy, P.; Sommers, E.: 
A new family of isolated symplectic singularities with trivial local fundamental group, 
Proc. Lond. Math. Soc. (3) {\bf 126} (2023), no. 5, 1496 - 1521.

 
\bibitem[BD]{BD} Bielawski, R.,  Dancer, A.: 
The geometry and topology of toric hyperk\"{a}hler manifolds.
Comm. Anal. Geom. {\bf 8} (2000), no. 4, 727 - 760


\bibitem[B-F]{B-F} Bielawski, R., Foscolo, L.: 
Hypertoric varieties, W-Hilbert schemes, and Coulomb branches, arXiv: 2304.08125


 
\bibitem[Dr]{Dr} Druel, S.:
Structures de contact sur les vari\'{e}t\'{e}s alg\'{e}briques de dimension 5. C. R. Acad. Sci. Paris Ser. I Math. {\bf 327} (1998), no. 4, 365 - 368.


\bibitem[Go]{Go} Goto, R.: On toric hyper-K\"{a}hler manifolds given by the hyper-Kaehler quotient method, in Infinite Analysis, World Scientific, 1992, pp. 317 - 338


\bibitem[HL]{HL} Hamm, H. A.,  L\^{e} D\~{u}ng Tr\'{a}ng.:
Rectified homotopical depth and Grothendieck conjectures. The Grothendieck Festschrift, Vol. II, 311 - 351,
Progr. Math., {\bf 87}, Birkh\''{a}user Boston, Boston, MA, 1990.


\bibitem[HS]{HS} Hausel, T.; Sturmfels, B.: 
Toric hyperK\"{a}hler varieties. 
Doc. Math. {\bf 7} (2002), 495 - 534.

     
\bibitem[KPSW]{KPSW} Kebekus, S., Peternell, T., Sommese, A.J, Wi\'{s}niewski, J.A: 
Projective contact manifolds, 
Invent. Math. {\bf 142} (2000), no. 1, 1 - 15.


\bibitem[Ko]{Ko} Konno, H.: 
Cohomology rings of toric hyperkaehler manifolds.
Internat. J. Math. {\bf 11} (2000), no. 8, 1001 - 1026.


\bibitem[M]{M} Milnor, J.:
Singular points of complex hypersurfaces, 
Annals of Mathematics Studies, No. 61. Princeton University Press, Princeton, NJ; University of Tokyo Press, Tokyo, 
1968. iii+122 pp.


\bibitem[Na 1]{Na 1}  Namikawa, Y.: 
A finiteness theorem on symplectic singularities, 
Compos. Math. {\bf 152} (2016), no. 6, 1225 - 1236.


\bibitem[Na 2]{Na 2} Namikawa, Y.: 
A characterization of nilpotent orbit closures among symplectic singularities, 
Math. Ann. {\bf 370} (2018), no. 1-2, 811 - 818.

 
\bibitem[Na 3]{Na 3} Namikawa, Y.: 
Fundamental groups of symplectic singularities. Higher dimensional algebraic geometry -in honour of Professor Yujiro Kawamata's sixtieth birthday, 321 - 334,
Adv. Stud. Pure Math., {\bf 74}, Math. Soc. Japan, Tokyo, 2017.


\bibitem[Ye]{Ye} Ye, Y-G.:
A note on complex projective threefolds admitting holomorphic contact structures, 
Invent. Math. {\bf 115} (1994), no. 2, 311 - 314.

\end{thebibliography}
\end{document}